\documentclass{elsart}

\usepackage{latexsym}
\usepackage[psamsfonts]{amssymb}
\usepackage{amsfonts}
\usepackage{graphicx,amssymb,enumerate}
\usepackage{amsmath}

\usepackage{subfig}
 \usepackage[misc]{ifsym}

\usepackage{balance}

\newtheorem{tm}{Theorem}[section]   % number due to section
\newtheorem{pp}[tm]{Proposition}   %  [tm] in the middle, let all the def,remark, thm, use one sequence of number

\newtheorem{re}[tm]{Remark}

\begin{document}

\allowdisplaybreaks

\newcommand{\al}{\alpha}
\newcommand{\be}{\beta}
\newcommand{\var}{\varepsilon}
\newcommand{\la}{\lambda}
\newcommand{\de}{\delta}
\newcommand{\st}{\stackrel}
\newcommand{\f}{F^{*}(R)}
\newcommand{\ce}{\centerline}
\newcommand{\Llr}{\Longleftrightarrow}
\newcommand{\Lr}{\Longrightarrow}

\begin{frontmatter}



\title{Approximation capability of the convolution methods for fuzzy
numbers}
\author{Huan Huang}
\author{}\ead{hhuangjy@126.com }

\address{Department of Mathematics, Jimei
University, Xiamen 361021, China}

\date{}
\begin{abstract}
This paper shows that how to approximate general fuzzy number by using convolution method.
\end{abstract}


\end{frontmatter}

%\title{}

%\author{Qing-song Mao\\
%Jimei University\\ College of Teacher Education \\
%Xiamen,
%361021, China\\ pinem@163.com\\% For a paper whose authors are all at the same institution,
%\and
%First line of institution2 address\\ Second line of institution2 address\\
%SecondAuthor@institution2.com\\
%}

%\maketitle %\thispagestyle{empty}

%-------------------------------------------------------------------------

%-------------------------------------------------------------------------
\section{Instructions}

  The approximation of fuzzy numbers attract many peoples attentions \cite{gp,ro,yo,yo2}.
Since
differentiable fuzzy numbers play an important role in the
implementation of fuzzy intelligent systems and their applications
(see \cite{co,fs}).
 Chalco-Cano et al. \cite{yo, yo2} introduced a
method based on the convolution to approximate
non-differentiable fuzzy numbers by differentiable fuzzy
numbers.

A significant advantage of this method
is
that
it can generate differentiable fuzzy number such that the distance   between which
and the original fuzzy number is less than or equal to
arbitrary predetermine positive number.

However, in the previous work, not all the fuzzy numbers
can be approximated by this method.
So it is  natural for us to consider
the question
that
 can we use the convolution method to approximate general fuzzy numbers?

In this paper, we want to answer this question.

\section{Preliminaries}

\subsection{Fuzzy numbers}

Let $\mathbb{N}$ be the set of all natural numbers, $\Bbb R$ be the
set of all real numbers.  For details, we refer the reader to
references \cite{da,wu}.

A fuzzy subsets $u$ on $\Bbb R$ can be seen as a mapping  from
$\mathbb{R}$ to [0,1]. For $\al\in (0,1]$, let $[u]_{\al}$ denote
the $\al$-cut of $u$; i.e., \ $[u]_{\al}\equiv\{x\in \mathbb
R:u(x)\geq \al \}$ and $[u]_0$ denotes $\overline{\{x\in \mathbb
R:u(x)>0\}}$. We call $u$ a fuzzy number if $u$ has the following
properties:
\\
(\romannumeral1) $[u]_1\neq \emptyset$; and
\\
(\romannumeral2) $[u]_\al=[u^{-}(\al), u^{+}(\al)]$ are compact intervals of
$\mathbb{R}$ for all $\al\in [0,1]$.
\\
The set of all fuzzy numbers is denoted by
$\mathcal{F}(\mathbb{R})$. In \cite{yo}, a fuzzy number is also called a fuzzy number.

The following  is a widely used representation theorem of
 fuzzy numbers.
\begin{pp}
\label{ulc}(Goetschel and Voxman \cite{gn})
 Given $u\in \mathcal{F}(\mathbb{R}),$ then

(\romannumeral1) \ $u^{-}(\cdot)$ is a left-continuous nondecreasing bounded
function on $(0,1]$;\\
(\romannumeral2) \ $u^{+}(\cdot)$ is a left-continuous nonincreasing
bounded function on $(0,1]$;\\
(\romannumeral3) \ $u^{-}(\cdot)$ and  $u^{+}(\cdot)$ are right continuous at
$\al=0$;\\
(\romannumeral4) \ $u^{-}(1)\leq u^{+}(1).$

Moreover, if the pair of functions $a(\la)$ and $b(\la)$ satisfy
conditions (\romannumeral1) through (\romannumeral4), then there exists a unique $u\in \mathcal{F}(\mathbb{R})$ such
that $[u]_{\alpha}=[a(\la),b(\la)]$ for each $\alpha\in (0,1].$
\end{pp}

The algebraic operations
 on $\mathcal{F}(\mathbb{R})$ are defined  as follows: given $u,\
v\in \mathcal{F}(\mathbb{R})$, $\al\in [0,1]$,
  \begin{gather*}
[u+v]_\al= [u]_\al+[v]_\al=[u^-(\alpha) + v^-(\alpha),\ u^+(\alpha) + v^+(\alpha)],\nonumber
\\
[u-v]_\al= [u]_\al-[v]_\al=[u^-(\alpha) - v^+(\alpha),\ u^+(\alpha) - v^-(\alpha)], \nonumber
\\
[u\cdot v]_\al=[u]_\al\cdot[v]_\al
=
[
\min  \{    xy:     \;     x\in [u]_\al,\, y\in[v]_\al   \}
,\
 \max  \{     xy:    \;     x\in [u]_\al,\, y\in [v]_\al   \}
 ].
\end{gather*}

The supremum metric
on $\mathcal{F}(\mathbb{R})$
is defined by
$$
d_{\infty}(u,v)=\sup_{\al\in [0,1]} \max\{|u^{-}(\al)-v^{-}(\al)|,\
|u^{+}(\al)-v^{+}(\al)|\},
$$
where $u,v\in
\mathcal{F}(\mathbb{R}).$

\subsection{Convolution of fuzzy numbers}

In this paper, we want to discuss the properties of
 sup-min convolution $u\nabla v$ of fuzzy numbers $u$ and $v$, which is defined by
$$(u\nabla v)(x)=\sup_{y\in \mathbb{R}}\{u(y)\wedge v(x-y)\}.$$

\begin{re}

In fact $u\nabla v=u+v$ for all $u,v\in \mathcal{F}(\mathbb{R})$. For details, see \cite{da, wu}.

\end{re}

%
%
%
%%%%%%%%%%%%%%%%%%%%%%%%%%%%%%%%%
\iffalse
It is known that for each fuzzy number $u$, there are four real
numbers $a$, $b$, $c$ and $d$ and two upper semicontinuous functions
$L:[a, b]\to [0, 1]$ and $R:[c, d]\to [0,1]$ with $L(a)=R(d)=0$,
$L(b)=R(c)=1$, and $L$ is nondecreasing and $R$ is nonincreasing,
such that $u$ can be described in the following matter (\cite{dp}):
\begin{equation} \label{1}u(x)=\left\{
\begin{array}{ll}
0, & \mbox{if}\ x<a, \\
L(x),& \mbox{if}\ a\leq x< b,\\
1, & \mbox{if}\ b\leq x\leq c,\\
R(x),& \mbox{if}\ c\leq x\leq d,\\
0, & \mbox{if}\ d<x.
\end{array}
\right.
\end{equation}
$a=c$, $b=d$ are admissible. When $a=c$, for instance, then
\[ u(x)=\left\{
\begin{array}{ll}
R(x),& \mbox{if}\ a\leq x\leq d,\\
0, & \mbox{if}\ x<a \ \mbox{or}\ d<x.
\end{array}
\right.
\]

\fi
%%%%%%%%%%%%%%%%%%%%%%%%%%%%%%%%%%%
%
%
%

In the following, we list some symbols which are used to denote subsets of
$\mathcal{F}(\mathbb{R})$.
\begin{itemize}

\item \
$\mathcal{F}_\mathrm{T}(\mathbb{R}) $ is denoted
the family of all fuzzy numbers $u$ such that $u$ is strictly
increasing on
$[u^-(0), u^-(1)]$,
strictly decreasing on
$[u^+(1), u^+(0)]$,
 and
  differentiable on
   $(  u^-(0), u^-(1)   ) \cup  (   u^+(1), u^+(0)   )$.

\vspace{1mm}

\item \
$\mathcal{F}_\mathrm{N}(\mathbb{R})$ is denoted the family of all
fuzzy numbers $u$ such that $u$ is differentiable on $(u^-(0), u^-(1)) \cup  (u^+(1), u^+(0))$.

\vspace{1mm}

\item \
$\mathcal{F}_\mathrm{C}(\mathbb{R})$ is denoted the family of all
continuous fuzzy numbers, i.e., the family of all fuzzy
numbers $u$ such that $u:\mathbb{R} \to [0, 1]$ is continuous
on $(u^-(0), u^+(0))$.

\vspace{1mm}

\item \
$\mathcal{F}_\mathrm{D}(\mathbb{R})$ is denoted the family of all
differentiable fuzzy numbers, i.e., the family of all fuzzy
numbers $u$ such that $u:\mathbb{R} \to [0, 1]$ is differentiable
on $(u^-(0), u^+(0))$.

\end{itemize}
Given
 a
 fuzzy number $u$ in $\mathcal{F}_\mathrm{N}(\mathbb{R})$,
$u$ need not be strictly increasing on $(u^-(0), u^-(1))$  and
strictly decreasing on $(u^+(1), u^+(0))$. So
$$\mathcal{F}_\mathrm{T}(\mathbb{R})
 \subsetneq
\mathcal{F}_\mathrm{N}(\mathbb{R}).$$

Observe
 that
  $p$ is differentiable on $(p^-(1),   p^+(1))$ for all $p\in \mathcal{F}(\mathbb{R})$.
Thus
we
know that, for each $u\in \mathcal{F}_\mathrm{N}(\mathbb{R})$,
 the   possible non-differentiable points are $u^-(1)$ and $u^+(1)$.

It
is easy to check that
\begin{gather*}
 \mathcal{F}_\mathrm{D}(\mathbb{R})
 \subsetneq
  \mathcal{F}_\mathrm{C}(\mathbb{R})
 \cap
      \mathcal{F}_\mathrm{N}(\mathbb{R}).
\end{gather*}

 Chalco-Cano et al. \cite{yo} constructed
fuzzy numbers $w_p,\ p>0$, which are defined by
\begin{equation} \label{wpc}
 w_p(x)=\left\{
\begin{array}{ll}
1-{\left(\frac{x}{p}\right)}^2,& \mbox{if}\ x\in [-p,p],\\
0, & \mbox{if}\ x\notin [-p,p].
\end{array}
\right.
\end{equation}
Obviously, $w_p\in \mathcal{F}_\mathrm{D}(\mathbb{R})$ for all  $p>0$.
They
 presented the
 following
 result.
\begin{pp} \cite{yo} \label{1}
If
$u\in \mathcal{F}_\mathrm{T}(\mathbb{R})$,
then
$u\nabla w_p \in \mathcal{F}_\mathrm{D}(\mathbb{R})$.
\end{pp}
Notice
that
 $d_\infty (  u, u\nabla w_{p}   )   \to    0$
as $p\to 0$.
Thus Proposition \ref{1}
indicates
that
every fuzzy number in $\mathcal{F}_\mathrm{T}(\mathbb{R})$ can be approximated by
fuzzy numbers sequences in $\mathcal{F}_{\mathrm{D}}(\mathbb{R})$.

We can see that
the fuzzy numbers $w_p$,  $p>0$, work as  smoothers, which transfer each fuzzy number $u$
to a differentiable (smooth) fuzzy number $u \nabla w_p$. This sequence of smooth fuzzy numbers construct a approximation
sequence of the original fuzzy number $u$,
i.e.,
 $u\nabla w_{p}    \to  u$ as $p\to 0$.

 Chalco-Cano et al. \cite{yo2}
  further putted forward a method to define smoothers.
  Suppose
that $p>0$ is a real number
 and
     that   $f: [0,1]   \to  [0,1]$ is a continuous and strictly decreasing function with $f(0)=1, f(1)=0$.
      A class of fuzzy numbers $Z_p^f$ is defined by
\[
Z_p^f(x)=
\begin{cases}
f^{-1} ( \|x\| / p), &   \|x\|  \leq   p,
\\
0,      &     \|x\|  >   p.
\end{cases}
\]
It is easy to show that $Z_p^f=w_p$ when $f=\sqrt{1-t}$.
They gave the following
result.
\begin{pp} \cite{yo2}
If $f$ is differentiable and $\lim_{\alpha \to 1-} f'(\al)   =  -\infty$,
then
$u \nabla Z_p^f    \in   \mathcal{F}_\mathrm{D}(\mathbb{R})$
for each $u\in   \mathcal{F}_\mathrm{T}(\mathbb{R})$.
\end{pp}
Notice
that $d_\infty ( u \nabla Z_p^f, \   u  ) \to  0 $ as $p\to 0$.
This means that
given $f$ satisfies the above conditions,
we obtain a smooth approximation
$\{  u \nabla Z_p^f:  \  p>0 \}$ of the fuzzy number $u$.
Different $f$ corresponds to different sequence of smooth approximation.

\section{Approximation}

In the previous work, only fuzzy numbers in $\mathcal{F}_\mathrm{T}(\mathbb{R})$
can be approximated. This type of fuzzy numbers have at most two possible non-differentiable points:
the endpoints of the 1-cut. Whereas, an arbitrarily given fuzzy number may have other non-differentiable points
or non-continuous points.
So it is  natural for us to consider
the question
that
 can we use the convolution method to approximate general fuzzy numbers?

In this section, to answer this question,
it discusses how to choose smoothers to smooth an arbitrarily given fuzzy number.
The following theorems are devoted to this problem.

\begin{tm}
\label{rap}
Suppose that $u\in \mathcal{F}_\mathrm{N}(\mathbb{R})\cap \mathcal{F}_\mathrm{C}(\mathbb{R})$   and that
 $w\in \mathcal{F}_\mathrm{D}(\mathbb{R})$, then $u\nabla
w\in \mathcal{F}_\mathrm{D}(\mathbb{R})$ when $w$ satisfies the following conditions.
\begin{enumerate}[(i)]
  \item $w(w^-(0))=u(u^-(0))$ and $w(w^+(0))=u(u^+(0))$.

    \item If $u^-(1)$ is the inner point of $[u]_0$, then
$w_-'(w^-(1))=0$.

  \item If $u^+(1)$ is the inner point of $[u]_0$, then
$w_+'(w^+(1))=0$.
\end{enumerate}

\end{tm}

\begin{tm}
\label{rndp}
Suppose that $u\in \mathcal{F}_\mathrm{C}(\mathbb{R})$
and
that $w\in      \mathcal{F}_\mathrm{D}(\mathbb{R})$,
then
$u\nabla w   \in \mathcal{F}_\mathrm{D}(\mathbb{R})$
when
$w$ satisfies the conditions (\romannumeral1), (\romannumeral2), and (\romannumeral3) in  Theorem \ref{rap}
and
the following condition (\romannumeral4).

(\romannumeral4) \
Given  $ x  \in  (u^-(0), u^+(0)) $ with $u(x)=\al<1$,
if
$x$
 is a non-differentiable point of $u$,
 then
  \\
(\romannumeral4) (1) \
if $x< u^-(1)$, then $w'(w^-(\al))=0$;
\\
(\romannumeral4) (2) \
if $x> u^+(1)$, then $w'(w^+(\al))=0$.

\end{tm}

\begin{tm}
\label{rncp}
Suppose that $u\in \mathcal{F}(\mathbb{R})$
and
that $w\in      \mathcal{F}_\mathrm{D}(\mathbb{R})$,
then
$u\nabla w   \in \mathcal{F}_\mathrm{D}(\mathbb{R})$
when
$w$ satisfies the conditions (\romannumeral1), (\romannumeral2), and (\romannumeral3) in  Theorem \ref{rap},
the condition
 (\romannumeral4) in Theorem \ref{rndp},
and
the following condition (\romannumeral5).

(\romannumeral5) \
Given  $ x  \in  (u^-(0), u^+(0)) $ with $u(x)=\al<1$,
if
$x$
 is a non-continuous point of $u$,
 then
  \\
(\romannumeral5) (1) \
if $x< u^-(1)$, then $w'(w^-(\beta))=0$, where $\beta=\lim_{y\to x-}   u(y)$;
\\
(\romannumeral5) (2) \
if $x> u^+(1)$, then $w'(w^+(\gamma))=0$,  where $\gamma=\lim_{z\to x+}   u(z)$.

\end{tm}

By above theorems, given an arbitrary fuzzy number,
if
 the number of non-continuous points of $u$ in
$[u]_0$
is finite,
then
we can construct a sequence of smoothers $v_n$, $n\in \mathbb{N}$,
such that
$u\nabla v_n $ is a differentiable fuzzy number for all   $n$,
and that
$\{u\nabla v_n \}$
converges to
$u$ in the supremum metric $d_\infty$.

\end{document}